\documentclass[10pt]{report}
\usepackage{amsmath}
\usepackage[cp866]{inputenc}
\usepackage[english,russian]{babel}
\usepackage{amsfonts,amssymb}
\usepackage{latexsym}
\usepackage{graphicx}

\def\qed{\hfill $\squar$}
\def\squar{\vbox{\hrule\hbox{\vrule height 6pt \hskip 6pt\vrule}\hrule}}

\newcommand{\eps}{\varepsilon}

\usepackage[english,russian]{babel} 
\usepackage{amssymb,amsmath} 

\newtheorem{theorem}{Теорема}[chapter]
\newtheorem{lemma}{Лемма}[chapter]

\newtheorem{corollary}{Следствие}[chapter]

\begin{document} \large
\setcounter{equation}{0}
\thispagestyle{empty}


\setcounter{secnumdepth}{4} \setcounter{tocdepth}{2}



\begin{center}{\bf \Large УСТОЙЧИВОСТЬ НЕЗАЛИПАЮЩИХ ПЕРИОДИЧЕСКИХ КОЛЕБАНИЙ, УСТАНОВЛЕННЫХ
МЕТОДОМ УСРЕДНЕНИЯ В РАЗРЫВНЫХ СИСТЕМАХ. I. ГЛАДКИЕ ВНЕ РАЗРЫВОВ
СИСТЕМЫ}\\
О.Ю. Макаренков
\end{center}

{\bf Аннотация.} {\small В настоящей работе утверждение второй
теоремы Н.Н. Боголюбова о периодических решениях гладких систем с
малым параметром обосновывается для разрывных систем, в которых
порождающее решение пересекает гиперплоскости разрыва
трансверсально и которые непрерывно дифференцируемы вне этих
гиперплоскостей. Данная ситуация имеет место в системах с сухим
трением в отсутствии залипания и упругих ограничителей. В качестве
иллюстрации доказывается устойчивость колебаний скорости тела,
перемещаемого вибрациями.}

\section{Введение}

Рассмотрим систему
\begin{equation}\label{ps}
  \dot x+h(t,x)=\eps f(t,x,\eps),
\end{equation}
где $h\in C^1(\mathbb{R}\times\mathbb{R}^n,\mathbb{R}^n)$ и $f$ --
$T$-периодическая по времени непрерывно дифференцируемая функция,
терпящая разрывы 1-го рода в таких точках $(t,x,\eps),$ в которых
некоторые компоненты $x$ обращаются в нуль (см. условие (A1)
ниже). В предположении, что порождающая система
\begin{equation}\label{np}
  \dot x+h(t,x)=0
\end{equation}
допускает только $T$-периодические решения, настоящая статья
изучает существование, единственность и устойчивость
$T$-периодических решений системы (\ref{ps}). Известной
$T$-периодической заменой переменных (см. замену \ref{ZAM} ниже)
система (\ref{ps}) приводится к стандартной форме принципа
усреднения. Соответственно, в случае, когда $f\in
C^1(\mathbb{R}\times\mathbb{R}^n\times[0,1],\mathbb{R}^n),$
поставленная задача полностью решена Н.Н. Боголюбовым в его второй
теореме (см. \cite{bog}, Ч.1, \S5, Теорема~II).

Принцип усреднения для нахождения периодических колебаний в
разрывных системах до сих пор применялся либо без обоснования (см.
\cite{bolotnik}), либо без обоснования устойчивости (см.
\cite{fec}), либо на основании негладкого аналога первой теоремы
Н.Н.~Боголюбова (см. \cite{babitski}, \cite{fidlin}, \cite{thom}).
Такой аналог впервые предложен В.А.~Плотниковым \cite{plo} и
позволяет  убедиться, что динамика системы (\ref{ps}) близка к
$T$-периодической на временном интервале порядка $[0,1/\eps]$ и не
гарантирует, что динамика системы действительно $T$-периодическая
и, тем более, устойчивая на всем $[0,+\infty).$ Результат
настоящий статьи впервые гарантирует последнее свойство. В
качестве иллюстративного примера в работе доказывается
$T$-периодичность и устойчивость колебаний скорости тела,
перемещаемого под действием периодических вибраций, что было ранее
установлено А. Фидлиным \cite{fidlin} на интервале $[0,1/\eps].$

Предлагаемый результат получен прямым методом склейки оператора
сдвига по траекториям системы (\ref{ps}) из его фрагментов на
гладких областях. Используя решения вспомогательных задач Коши с
векторным временем (см. лемму~\ref{lem1}), доказано, что оператор
Пуанкаре системы (\ref{ps}) дифференцируем по фазовой переменной и
параметру $\eps>0.$ Далее, используя сходимость правой части при
уменьшении $\eps>0$ по мере (см. следствие~\ref{cor2}),
установлено равенство классической функции усреднения и
производной оператора Пуанкаре по $\eps$ в $\eps=0.$ Это позволило
связать свойства собственных значений нулей функции усреднения с
такими свойствами оператора Пуанкаре,  которые достаточны для
анализа устойчивости его неподвижных точек методами теории
динамических систем.

\section{Основной результат}
На протяжении статьи  $\xi^j$ является $j$-й компонентой вектора
$\xi\in\mathbb{R}^n,$ $x(\cdot,\xi,0)$ обозначает решение
порождающей системы (\ref{np}) с
    начальным условием $x(0)=\xi$ и $B_r(\zeta)$ -- это шар в
    $\mathbb{R}^n$ радиуса $r>0$ с центром в точке
    $\zeta\in\mathbb{R}^n.$
Результат статьи применим к разрывным системам, удовлетворяющим
следующим аналогичным предположениям А.~Фидлина \cite{fidlin}
условиям.

\begin{itemize}
  \item[(A1)] Положим $\mathbb{R}^n_s=\left\{\xi\in\mathbb{R}^n:{\rm sign}(\xi^j)=s^j,j\in
    \overline{1,n}\right\},$
$s\in\{-1,1\}^n=\underbrace{\{-1,1\}\times\ldots\times\{-1,1\}}_{n\
  \mbox{\footnotesize{штук}}}.$  Существует $2^n$ функций $f_s\in C^1(\mathbb{R}\times\mathbb{R}^n\times[0,1],
  \mathbb{R}^n),$ $s\in \{-1,1\}^n$ таких, что
  $$
    f(t,\xi,\eps)=f_s(t,\xi,\eps),\quad
    (t,\xi,\eps)\in\mathbb{R}\times\mathbb{R}^n_s\times[0,1],\
    s\in\{-1,1\}^n.
  $$
\end{itemize}

Следующие два условия предъявляются к такому $x(\cdot,\xi_0,0),$
которое, ожидается, будет порождающим, но часто они выполнены или
нет сразу для всех решений системы (\ref{np}).

\begin{itemize}
  \item[(A3)] Множество точек $S\subset\mathbb{R}^n\backslash\cup_{s\in\{-1,1\}^n}\mathbb{R}^n_s,$
  в которых функция $\xi\mapsto f(t,\xi,\eps)$ не
  является непрерывно дифференцируемой, не зависит от $t$ и
  $\eps,$ и для любых $j\in\overline{1,n}$ и  $t\in\mathbb{R}$ существует $p\in\overline{1,n}$
такое, что функция $\xi\mapsto f^j(t,\xi^1,...,\xi^{p-1},0\cdot
\xi^p,\xi^{p+1},...,\xi^n,\eps)$ непрерывно дифференцируема в
точке
  $x(t,\xi_0,0).$
  \item[(A2)] Предположим, что множество $\{t\in[0,T]:x(t,\xi_0,0)\}\in S$ конечно и занумеруем его
  элементы как $0\le t_1<...<t_m<T.$ Пусть $t_1>0$ и для любых
  $j\in\overline{1,n}$ и $i\in\overline{1,m}$ таких, что
  $x^j(t_i,\xi_0,0)=0$ и $\{\xi\in\mathbb{R}^n:\xi^j=0\}\subset S,$
  имеем $(x^j)'_t(t_i,\xi_0,0)\not=0.$
\end{itemize}

Поскольку система (\ref{ps}) может вообще не иметь
дифференцируемого на всем временном промежутке решения, нам
следует принять несколько более общее определение.

{\bf Определение 1.} {\it Решением системы (\ref{ps}) называется
непрерывная функция $x,$ дифференцируемая всюду, за исключением,
быть может, множества $\{t:x(t)\in S\}$ и удовлетворяющая всюду,
кроме, быть может, этого множества, системе (\ref{ps}).}

Данное определение позволяет не ограничивая общности считать, что
функция $f$ ограничена на каждом ограниченном множестве.

\begin{lemma} \label{lem2} Пусть $h\in
C^1(\mathbb{R}\times\mathbb{R}^n,\mathbb{R}^n)$ и $f$
удовлетворяет условию (A1).  Пусть $\xi_0\in\mathbb{R}^n$ таково,
что выполнены условия (A2)-(A3). Тогда существует $\delta>0$
такое, что при всех $\eps\in[0,\delta],$ $v\in B_\delta(\xi_0)$
система (\ref{ps}) имеет единственное решение $t\mapsto
x(t,\xi,\eps)$ с начальным условием $x(0,\xi,\eps)=\xi.$ Это
решение продолжимо на $[0,T]$ и непрерывно дифференцируемо по
$(\xi,\eps)\in B_\delta(\xi_0)\times [0,\delta).$ Кроме того,
$$
  \left\{t:x^j(t,\xi,\eps)=0\right\}\subset\bigcup_{i=1}^m\left\{T_i^j(\xi,\eps)\right\},\quad
  j\in\overline{1,n},
$$
где $T_i^j\in C^1(B_\delta(\xi_0)\times B_\delta(0),\mathbb{R}^n)$
и $T^j_i(\xi_0,0)=t_i$ при всех $i\in\overline{1,m}$ и
$j\in\overline{1,n}.$
\end{lemma}

Для доказательства леммы~\ref{lem2} нам понадобится следующее
вспомогательное утверждение, в котором $1_{\mathbb{R}^n}$ -- это
$n$-мерный вектор, состоящий из единиц, и для произвольных
$g:\mathbb{R}^n\to\mathbb{R}^n,$ $t\in\mathbb{R},$
$\xi\in\mathbb{R}^n$ запись
$g(t1_{\mathbb{R}^n}+\overrightarrow{\xi})$ обозначает следующее
$$
  g(t1_{\mathbb{R}^n}+\overrightarrow{\xi})=\left(\begin{array}{l}
  g^1(t+\xi^1)\\
  \vdots\\
  g^n(t+\xi^n)\end{array}\right).
$$

\begin{lemma} \label{lem1} Пусть $F\in
C^1(\mathbb{R}\times\mathbb{R}^n\times[0,1],\mathbb{R}^n).$ Тогда
для любых $t_*\in\mathbb{R}$ и $\xi_*\in\mathbb{R}^n$ существует
$\gamma>0$ такое, что при любых $\Delta\in B_\gamma(0),$ $\xi\in
B_\gamma(\xi_*),$ $\eps\in[0,\gamma),$ задача
\begin{eqnarray}
 & & \dot x=F(t,x,\eps),\label{S}\\
 & & x(t_*1_{\mathbb{R}^n}+\overrightarrow{\Delta})=\xi_*\label{S1}
\end{eqnarray}
имеет единственное решение $t\mapsto
x(t,t_*1_{\mathbb{R}^n}+\overrightarrow{\Delta},\xi,\eps),$
определенное на $\mathbb{R}.$ Более того, функция $x$ непрерывно
дифференцируема на $\mathbb{R}\times
B_\gamma(t_*1_{\mathbb{R}^n})\times
B_\gamma(\xi_*)\times[0,\gamma).$
\end{lemma}

{\bf Доказательство.} Обозначим через
$\widetilde{x}(\cdot,t_*,\zeta,\eps)$ решение системы (\ref{S}) с
начальным условием $x(t_*)=\zeta.$ Рассмотрим функцию
$$
  \Phi(\Delta,\zeta,\xi,\eps)=\widetilde{x}(t_*+\overrightarrow{\Delta},t_*,\zeta,\eps)-\xi,
$$
непрерывно дифференцируемую на
$\mathbb{R}^n\times\mathbb{R}^n\times\mathbb{R}^n\times[0,1].$
Имеем $\Phi(0,\xi_*,\xi_*,0)=0$ и ${\rm
det}\|\Phi'_\xi(0,\xi_*,\xi_*,0)\|={\rm det}\|I\|\not=0.$ Поэтому,
из теоремы о неявной функции следует существование $\gamma>0$ и
функции $(\Delta,\xi,\eps)\mapsto \zeta(\Delta,\xi,\eps),$
непрерывно дифференцируемой на $B_\gamma(0)\times
B_\gamma(\xi_*)\times[0,\gamma)$ и удовлетворяющей условию
$$
  \Phi(\Delta,\zeta(\Delta,\xi,\eps),\xi,\eps)=0,\quad
  \|\Delta\|<\gamma,\ \|\xi-\xi_*\|<\gamma,\ \eps\in[0,\gamma).
$$
Искомая функция $x$ дается формулой
$$
  x(t,t_*1_{\mathbb{R}^n}+\overrightarrow{\Delta},\xi,\eps)=\widetilde{x}(t,t_*,\zeta(\Delta,\xi,\eps),\eps),
$$
и справедливость соотношения (\ref{S1}) проверяется
непосредственно.\qed

{\bf Доказательство леммы~\ref{lem2}.} Доказательство проводится в
два этапа. На первом этапе строятся фрагменты решения системы
(\ref{ps}), проходящего в окрестности решения $x(\cdot,\xi_0,0)$ и
компоненты которого определены на подходящих окрестностях
интервалов $(t_{i-1},t_i).$ На втором этапе искомое решение
склеивается из полученных фрагментов.

{\bf I этап.} Обозначим $t_0=0,$ $t_{m+1}=T,$ зафиксируем
произвольное $i\in\overline{1,m+1},$ положим
$(t_*,\xi_*)=(t_{i-1},x(t_{i-1},\xi_0,0))$ и применим
лемму~\ref{lem1} к системе
\begin{equation}\label{SV}
 \dot x+h(t,x)=\eps f_s(t,x,\eps),
\end{equation}
где $s^j={\lim\limits_{t\to
 t_{i-1}+0}{\rm sign}\left(x^j(t,\xi_0,0)\right)},$
 $j\in\overline{1,n}.$ Пусть $\gamma>0$ -- то число, а
 $\widetilde{x}_i(t,t_{i-1}1_{\mathbb{R}^n}+\overrightarrow{\Delta},\xi,\eps)$ --
 та функция, о которых говорится в лемме~\ref{lem1}. Обозначим
 через $\nu_1,...,\nu_{k_i}$ номера компонент порождающего решения
 $x(\cdot,\xi_0,0),$ обращающиеся в нуль при $t=t_i.$ Из условия
 (A3) следует, что
 $$
   (\widetilde{x}_i^j)'_t(t_i,t_{i-1}1_{\mathbb{R}^n},\xi_*,0)=(x^j)'_t(t_i,\xi_0,0)\not=0,\quad
   j\in\{\nu_1,...,\nu_{k_i}\}.
 $$
 Поэтому, теорема о неявной функции позволяет утверждать, что
 существует $0<\delta_i<\gamma$ и $k_i$ функций
 $t_{i}^j:\mathbb{R}^n\times\mathbb{R}^n\times[0,1]\to\mathbb{R},$
 $j\in
\{\nu_1,...,\nu_{k_i}\},$
удовлетворяющих условиям (см. рис.~1):

a) $t_i^j\in C^1(B_\delta(t_{i-1}1_{\mathbb{R}^n})\times
B_\delta(\xi_*)\times[0,\delta),\mathbb{R}^n),$

б) $t_i^j(0,\xi_*,0)=t_i,$

в)
$\widetilde{x}_i^j(t_i^j(\Delta,\xi,\eps),t_{i-1}1_{\mathbb{R}^n}+\Delta,\xi,\eps)=0,$
 $\|\Delta\|<\delta_i,$ $\xi\in B_{\delta_i}(\xi_*),$ $\eps\in[0,\delta_i).$

\begin{figure}[htb]
\centerline{\includegraphics[scale=0.33]{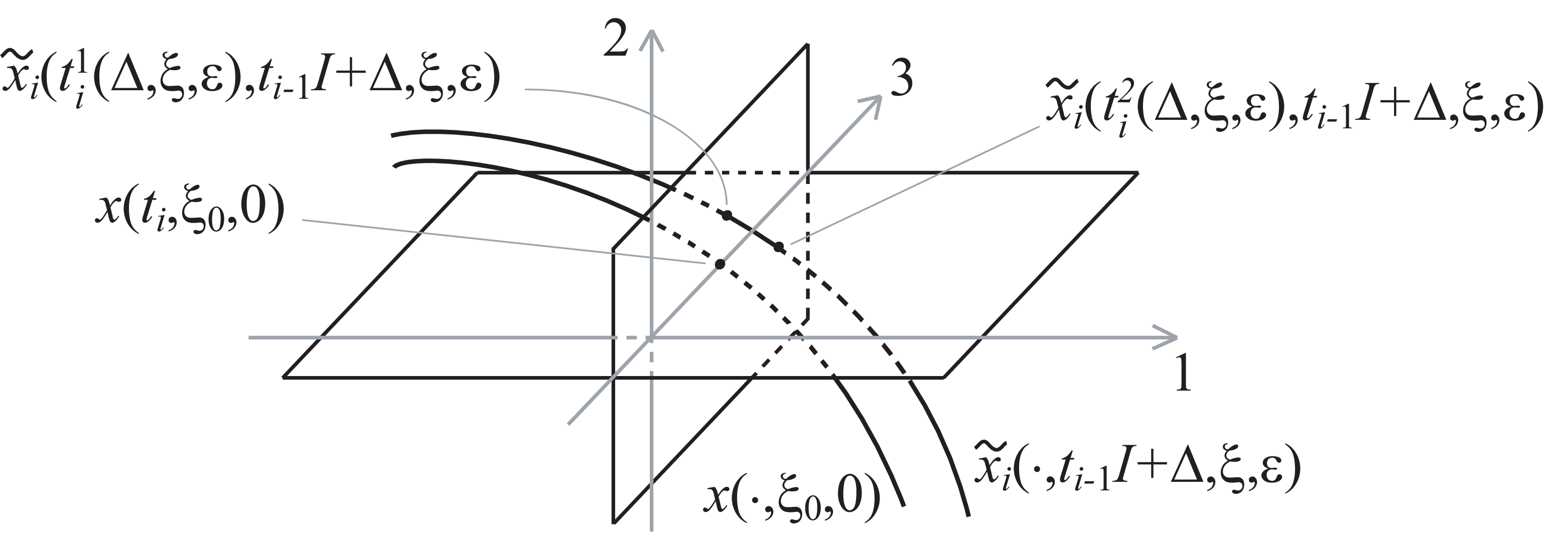}} 
{\footnotesize Рис. 1. Иллюстрация возможного поведения решения
$\widetilde{x},$ приводящего к разветвлению момента времени $t_i$
на две функции $t_i^1$ и $t_i^2,$ то есть к ситуации, когда разные
компоненты $\xi\mapsto f^{j_1}(t,\xi,\eps)$ и $\xi\mapsto
f^{j_2}(t,\xi,\eps)$ имеют разрывы в одной и той же точке
$x(t_i,\xi_0,0).$ Такой случай допускается условиями (A2)-(A3) и
реализуется в системах с трением с несколькими степенями свободы
(см.  \cite{bolotnik}).}\label{fig}
\end{figure}

В силу единственности неявной функции, $\delta_i>0$ может быть
уменьшено ещё и так, что
$\widetilde{x}^j_i(t,t_{i-1}1_{\mathbb{R}^n}+\Delta,\xi,\eps)\not=0$
при всех $t\in (t_{i-1}+\Delta^j,t_i^j(\Delta,\xi,\eps)),$
$j\in\{\nu_1,...,\nu_{k_i}\},$ $\|\Delta\|<\delta_i,$
$\|\xi-\xi_*\|<\delta_i,$ $0\le\eps<\delta_i.$ Для оставшихся
компонент $t_i$ определим как
\begin{equation}\label{OST}
t^j_i(\Delta,\xi,\eps):=t_i,\quad\Delta,\xi\in\mathbb{R}^n,\
\eps\in[0,1],\
j\in\overline{1,n}\backslash\{\nu_1,...,\nu_{k_i}\}\end{equation}
и уменьшим $\delta_i>0,$ если необходимо, так, что
$\widetilde{x}^j(t,t_{i-1}1_{\mathbb{R}^n}+\Delta,\xi,\eps)\not=0$
при всех $t\in(t_{i-1}+\Delta^j,t_i^j(\Delta,\xi,\eps)],$
$j\in\overline{1,n}\backslash\{\nu_1,...,\nu_{k_i}\},$
$\|\Delta\|<\delta_i,$ $\|\xi-\xi_*\|<\delta_i,$
$0\le\eps<\delta_i.$

{\bf II этап.} Двигаясь от $i=m+1$ до $i=2,$ уменьшим
$\delta_{i-1}>0$ одно за другим так, чтобы
$$
\begin{array}{l}
  \|\widetilde{x}_{i-1}(t_{i-1}(\Delta,\xi,\eps),t_{i-2}1_{\mathbb{R}^n}+\Delta,\xi,\eps)-x(t_{i-1},\xi_0,0)\|<\delta_i,\\
  \|t_{i-1}(\Delta,\xi,\eps)-t_{i-1}\|<\delta_i,\\
  \mbox{при всех} \|\Delta\|<\delta_{i-1},
 \|\xi-x(t_{i-2},\xi_0,0)\|<\delta_{i-1},
0\le\eps<\delta_{i-1}. \end{array}$$
  В помощь читателю мы подробно выписываем
первые итерации построения функции $x,$ но мелким шрифтом.
{\footnotesize Итак, для каждого $j\in\overline{1,n},$
$\|\xi-\xi_0\|\le\delta_1$ и $0\le\eps<\delta_1$ положим
$$
\begin{array}{l}
  x^j(t,\xi,\eps):=\widetilde{x}^j_1(t,0,\xi,\eps),\\
  \mbox{при \ всех\ }t\in[0,t_1^j(0,\xi,\eps)],
\end{array}
$$
$$
\begin{array}{l}
  x^j(t,\xi,\eps):=\widetilde{x}^j_2(t,t_1(0,\xi,\eps),x(\overrightarrow{t}_1(0,\xi,\eps),\xi,\eps),\eps),\\
  \mbox{при \ всех\
  }t\in[t_1^j(0,\xi,\eps),t_2^j(t_1(0,\xi,\eps)-t_1,x(t_1(0,\xi,\eps),\xi,\eps),\eps)].
\end{array}
$$
Далее, используя обозначение
$$
  T_2^j(v,\eps)=t_2^j(t_1(0,\xi,\eps)-t_1,x(t_1(0,\xi,\eps),\xi,\eps),\eps),
$$
построение продолжается как
$$
\begin{array}{l}
  x^j(t,\xi,\eps):=\widetilde{x}^j_3(t,T_2(\xi,\eps),x(\overrightarrow{T}_2(\xi,\eps),\xi,\eps),\eps),\\
  \mbox{при \ всех\
  }t\in[T_2^j(\xi,\eps),t_3^j(T_2^j(\xi,\eps)-t_2,x(T_2(\xi,\eps),\xi,\eps),\eps)].
\end{array}
$$}
\noindent Общая итерационная формула для определения
$x(t,\xi,\eps)$ при произвольных $t\in[0,T],$
$\|\xi-\xi_0\|<\delta_1,$ $0\le\eps<\delta_1$ и $i=1,...,m+1$
выписывается как
$$
\begin{array}{l}
  x^j(t,\xi,\eps):=\widetilde{x}^j_i(t,T_{i-1}(\xi,\eps),x(\overrightarrow{T}_{i-1}(\xi,\eps),\xi,\eps),\eps),\\
  \mbox{при \ всех\
  }t\in[T_{i-1}^j(\xi,\eps),T_i^j(\xi,\eps)],
\end{array}
$$
где
$T_i(\xi,\eps)=t_i(T_{i-1}(\xi,\eps)-t_{i-1},x(T_{i-1}(\xi,\eps),\xi,\eps)),$
$T_1(\xi,\eps)=t_1(0,\xi,\eps),$ $T_0(\xi,\eps)=0.$ При этом из
(\ref{OST}) имеем $T_{m+1}^j=T$ для любого $j\in\overline{1,n}.$
Так как $\widetilde{x}_i\in C^1(\mathbb{R}\times
B_{\delta_1}(t_{i-1}1_{\mathbb{R}^n})\times
B_{\delta_1}(x(t_{i-1},\xi_0,0))\times[0,\delta_1),\mathbb{R}^n)$
и $t_i\in C^1(B_{\delta_1}(0)\times
B_{\delta_1}(x(t_{i-1},\xi_0,0))\times[0,\delta_1),\mathbb{R}^n),$
то $z(t,\cdot),T_i\in
C^1(B_{\delta_1}(x(t_{i-1},\xi_0,0))\times[0,\delta_1),\mathbb{R}^n)$
при всех $i\in\overline{1,m},$ $t\in[0,T].$

Для завершения доказательства нам остается обосновать
единственность построенного решения. Для этого достаточно
показать, что при любом $i\in\overline{1,m}$ и достаточно малом
$\gamma>0$ часть $x((t_*,t_*+\gamma))$ решения $x$ системы
(\ref{ps}) с начальным условием $x(t_*)=\xi,$ где
$|t_*-t_i|<\gamma,$ $\|\xi-\xi_0\|<\gamma$ и
$\xi,x(t_i,\xi_0,0)\in S,$ лежит в том же множестве
$\mathbb{R}^n_s,$ что и $x((t_i,t_i+\gamma),\xi_*,0)$ (как это
имеет место для построенного решения $t\mapsto x(t,\xi,\eps)$). В
силу принятого определения решения системы (\ref{ps}) можем
считать, что решение $x$ непрерывно дифференцируемо на
$(t_*,t_*+\gamma).$ Но тогда, считая $\eps>0$ достаточно малым,
получаем, что значения $x'(t),$ $t\in (t_*,t_*+\gamma)$ и
$x'_t(t,\xi_*,0),$ $t\in (t_i,t_i+\gamma)$ сколь угодно близки.
Требуемое утверждение теперь легко следует из трансверсальности
$t\mapsto x(t,\xi_0,0)$ по отношению к $S$ в точке $t_i,$
вытекающей из $(A3).$ \qed

Лемма~\ref{lem2} позволяет ввести при малых $\eps>0$ и
$\xi\in\mathbb{R}^n$ близких к $\xi_0$ следующую функцию
\begin{equation}\label{ZAM}
   u(t,\xi,\eps)=x^{-1}(t,x(t,\xi,\eps),0),
\end{equation}
где $x^{-1}(t,\cdot,0)$ -- обратный к $x(t,\cdot,0)$ оператор (то
есть $x(t,x^{-1}(t,\xi,0),0)=x^{-1}(t,x(t,\xi,0),0)=\xi),$
существующий в силу гладкости порождающей системы (\ref{np}).
Замена (\ref{ZAM}) приводит (\ref{ps}) к стандартной форме
принципа усреднения
\begin{equation}\label{pss}
  \dot u=\eps(x'_u(t,u,0))^{-1}f(t,x(t,u,0),\eps).
\end{equation}
Решения системы (\ref{pss}) будем понимать в смысле определения~1.
В частности функция $t\mapsto u(t,\xi,\eps)$ является решением
системы (\ref{pss}) и, в силу леммы~\ref{lem2}, непрерывно
дифференцируемо по $(\xi,\eps)$ достаточно близким к $(\xi_0,0).$
Нам понадобится ряд свойств правой части системы (\ref{pss}),
которые мы сейчас выведем из леммы~\ref{lem2}.
\begin{corollary} \label{cor1} В условиях леммы~\ref{lem2} функция
 $$
   t\mapsto(x'_u(t,\xi,0))^{-1}f(t,x(t,\xi,0),0)
 $$
 суммируема на $[0,T]$ при всех $\|\xi-\xi_0\|<\delta.$
\end{corollary}
{\bf Доказательство.} Утверждение следует из суммируемости функции
$t\mapsto f(t,x(t,\xi,0),0),$ которая, в силу леммы~\ref{lem2},
непрерывна на $[0,T]$ всюду, кроме, быть может, точек
$\cup_{i\in\overline{1,m},j\in\overline{1,n}}\{T_i^j(\xi,0)\}.$
\qed

Данное следствие позволяет ввести в рассмотрение классическую
функцию усреднения
$$
  \overline{f}(\xi)=\int_0^T(x'_u(\tau,\xi,0))^{-1}f(\tau,x(\tau,\xi,0),0)d\tau,\quad
  \|\xi-\xi_0\|<\delta.
$$

\begin{corollary}\label{cor2} В условиях леммы~\ref{lem2} при всех $\xi\in
B_\delta(\xi_0)$ и $\sigma>0$ справедливо соотношение
$$
  \lim\limits_{\eps\to 0}{\rm
  mes}\left\{t\in[0,T]:\|f(t,x(t,u(t,\xi,\eps),0),\eps)-f(t,x(t,\xi,0),0)\|\ge\sigma\right\}=0.
$$
\end{corollary}

{\bf Доказательство.} Зафиксируем $j\in\overline{1,n},$
$\sigma>0,$ $\xi\in B_\delta(\xi_0)$ и $\gamma>0.$ Выберем
$\eps_0>0$ настолько малым, что
$$
  \|T_{i}^j(\xi,\eps)-T_{i}^j(\xi,0)\|<\gamma,\quad
  \eps\in[0,\eps_0],\ i\in\overline{1,m}.
$$
Обозначая $T_{0}^j(\xi,\eps)\equiv 0$ и $T_{m+1}^j(\xi,\eps)\equiv
T,$ при любом $i\in\overline{1,m+1}$
 имеем
$$
  f(t,x(t,u(t,\xi,\eps),0),\eps)\to f(t,x(t,\xi,0),0)\quad\mbox{при}\
  \eps\to 0
$$
равномерно на $[T_{i-1}^j(\xi,0)+\gamma,T_{i}^j(\xi,0)-\gamma],$ в
частности мы можем уменьшить $\eps_0>0$ настолько, что
$$
  \|f(t,x(t,u(t,\xi,\eps),0),\eps)-f(t,x(t,\xi,0),0)\|<\sigma,
$$
при всех $[T_{i-1}^j(\xi,0)+\gamma,T_{i}^j(\xi,0)-\gamma],$
$\xi\in B_\delta(\xi_0),$ $\eps\in[0,\eps_0],$
$i\in\overline{1,m+1}.$ Таким образом,
\begin{eqnarray*}
&& {\rm
  mes}\left\{t\in[0,T]:\|f(t,x(t,u(t,\xi,\eps),0),\eps)-f(t,x(t,\xi,0),0)\|\ge\sigma\right\}\le\\
&& \le  (m+1)\cdot 2\gamma,\quad\mbox{при всех}\  \xi\in
B_\delta(\xi_0),\ \eps\in[0,\eps_0].
\end{eqnarray*}
 Поскольку
$\gamma>0$ было выбрано произвольно, следствие доказано. \qed

\begin{corollary}\label{cor3} В условиях леммы~\ref{lem2} имеем
$\overline{f}(\xi)=x'_\eps(T,\xi,0),$ в частности функция
$\overline{f}$ непрерывно дифференцируема на $B_\delta(v_0).$
\end{corollary}

{\bf Доказательство.} Имеем
$$
  u(T,\xi,\eps)=\xi+\eps\int_0^T(x'_u(\tau,u(\tau,\xi,\eps),0))^{-1}f(\tau,x(\tau,u(\tau,\xi,\eps),0),\eps)d\tau,
$$
поэтому,
\begin{eqnarray*}
  u'_\eps(T,\xi,0)&=&\lim\limits_{\eps\to
 0}\frac{u(T,\xi,\eps)-u(T,\xi,0)}{\eps}=\\
 &=&\lim\limits_{\eps\to
 0}\int_0^T(x'_u(\tau,u(\tau,\xi,\eps),0))^{-1}f(\tau,x(\tau,u(\tau,\xi,\eps),0),\eps)d\tau.
\end{eqnarray*}
В силу непрерывности функции $u$ и ограниченности функции $f,$
подынтегральное выражение равномерно ограничено по $\tau\in[0,T],$
$\xi\in B_\delta(\xi_0)$ и $\eps\in[0,\delta).$ Значит,
следствие~\ref{cor2} позволяет применить теорему Лебега о
предельном переходе под знаком интеграла и прийти к заключению
$$
  \lim\limits_{\eps\to
  0}\int_0^T(x'_u(\tau,u(\tau,\xi,\eps),0))^{-1}f(\tau,x(\tau,u(\tau,\xi,\eps),0),\eps)d\tau=\overline{f}(\xi),
$$
завершающему доказательство. \qed

\begin{theorem}\label{th1} Пусть $h\in C^1(\mathbb{R}\times\mathbb{R}^n,\mathbb{R}^n)$ и каждое решение
порождающей системы (\ref{np}) $T$-периодично. Пусть $f$ --
$T$-периодическая по времени непрерывно дифференцируемая функция,
терпящая разрывы 1-го рода на $S,$  точнее, пусть  выполнено
условие (A1). Зададимся $\xi_0\in\mathbb{R}^n,$ удовлетворяющим
(A2), то есть таким, что при каждом $t\in[0,T]$ и
$j\in\overline{1,n}$ решение $t\mapsto x(t,\xi_0,0)$ порождающей
системы (\ref{np}) пересекает не более одной гиперплоскости
разрыва функции $f^j$ и такие пересечения происходят только при
$t\in(0,T).$ Пусть, наконец, решение $x(\cdot,\xi_0,0)$ пересекает
$S$ трансверсально, то есть выполнено условие (A3). Тогда имеют
место следующие утверждения:
\begin{itemize}
\item[1)] Если $\overline{f}(\xi_0)=0$ и ${\rm
det}\|\overline{f}'(\xi_0)\|\not=0,$ то существуют $\delta>0$ и
$\eps_0>0$ такие, что при $\eps\in(0,\eps_0)$ система (\ref{pss})
имеет единственное $T$-периодическое решение $u_\eps$ с начальным
условием $u_\eps(0)\in B_\delta(\xi_0).$ \item[2)] Если в условиях
пункта 1) все собственные значения матрицы $\overline{f}'(\xi_0)$
имеют отрицательные вещественные части, то решения
$\{u_\eps\}_{\eps\in(0,\eps_0)}$ асимптотически устойчивы.
\item[3)] Если в условиях пункта 1) хотя бы одно собственное
значение матрицы $\overline{f}'(\xi_0)$ имеет положительную
вещественную часть, то решения $\{u_\eps\}_{\eps\in(0,\eps_0)}$
неустойчивы.
\end{itemize}
\end{theorem}

{\bf Доказательство.} Положим
$$
  \overline{f}_\eps(\xi)=\int_0^T(x'_u(\tau,u(\tau,\xi,\eps),0))^{-1}f(\tau,x(\tau,u(\tau,\xi,\eps),0),\eps)d\tau,
$$
тогда
\begin{equation}\label{TE}
  u(T,\xi,\eps)=\xi+\eps\overline{f}_\eps(\xi)=u'_\xi(T,\xi,0)+\eps\overline{f}_\eps(\xi).
\end{equation}
Поэтому,
$$
   \frac{u'_\xi(T,\xi,\eps)-u'_\xi(T,\xi,0)}{\eps}=(\overline{f}_\eps)'(\xi).
$$
В силу леммы~\ref{lem2} имеем
$\dfrac{u'_\xi(T,\xi,\eps)-u'_\xi(T,\xi,0)}{\eps}\to
{{u}'_\xi}{}'_\eps(T,\xi,0)$ при $\eps\to 0$ равномерно по $\xi\in
B_\delta(\xi_0).$ Следовательно, учитывая заключение
следствия~\ref{cor3},
$$
  (\overline{f}_\eps)'(\xi)\to(\overline{f})'(\xi)\quad\mbox{при}\
  \eps\to 0
$$
равномерно по $\xi\in B_\delta(\xi_0).$

1) Начнем с доказательства утверждения 1). Другими словами,
требуется показать, что существует $\eps_0>0$ и $\delta>0$ такие,
что при $\eps\in(0,\eps_0)$ функция $\xi\mapsto u(T,\xi,\eps)-\xi$
имеет единственный нуль в $B_\delta(\xi_0).$ В силу формулы
(\ref{TE}) достаточно установить данное утверждение для функции
$\overline{f}_\eps(\xi).$ Но в силу условия 2) теоремы это
утверждение немедленно следует из теоремы о неявной функции.

2) Перейдем к вопросу об устойчивости найденных решений. Для этого
изучим собственные значения матрицы $u'_\xi(T,\xi_\eps,\eps).$
Имеем
$$
  u'_\xi(T,\xi_\eps,\eps)=I+\eps(\overline{f}_\eps)'(\xi_\eps).
$$
Предположим, что вещественные части всех собственных значений
матрицы $(\overline{f})'(\xi_0)$ отрицательны. Пусть $\lambda_0$
-- какое-нибудь собственное значение матрицы
$(\overline{f})'(\xi_0)$ и $\lambda_\eps$ -- какое-нибудь
собственное значение матрицы $(\overline{f})'(\xi_\eps),$
сходящееся при $\eps\to 0$ к $\lambda_0.$ Тогда
$1+\eps\lambda_\eps$ будет являться собственным значением матрицы
$I+\eps(\overline{f}_\eps)'(\xi_\eps).$ Но
$\lambda_\eps=\lambda_0+\delta_\eps,$ где $\delta_\eps\to 0$ при
$\eps\to 0,$ значит
$1+\eps\lambda_\eps=1+\eps\lambda_0+\eps\delta_\eps.$ Так как
${\rm Re}(\lambda_0)<0,$ то существует $\eps_0>0$ такое, что ${\rm
Re}(1+\eps\lambda_\eps)<0$ при $\eps\in(0,\eps_0].$ Таким образом,
при $\eps\in(0,\eps_0]$ собственные значения матрицы
$u'_\xi(T,\xi_\eps,\eps)$ лежат  в единичном шаре. Зафиксируем
$\eps\in[0,\eps_0]$ и обозначим через $\|\cdot\|_0$ такую норму в
$\mathbb{R}^n,$ что
$$\sup_{\|\zeta\|_0\le 1}\|u'_\xi(T,\xi_\eps,\eps)\zeta\|_0\le q<1$$ (см.
\cite{polo}, с.~90, лемма~2.2). Тогда найдется $\delta>0$ такое,
что
$$\sup_{\|\zeta\|_0\le
1}\|u'_\xi(T,\xi,\eps)\zeta\|_0\le \widetilde{q}<1\quad\mbox{для
всех}\ \xi\in B_\delta(\xi_\eps).
$$
Следовательно, будем иметь
$$
  \|u(T,\xi_1,\eps)-u(T,\xi_2,\eps)\|_0\le\widetilde{q}\|\xi_1-\xi_2\|_0,\quad
  \xi_1,\xi_2\in B_\delta(\xi_0),
$$
что означает (см. \cite{kraope}, лемма~9.2) асимптотическую
устойчивость периодического решения $u_\eps.$

3) Пусть теперь матрица $(\overline{f})'(\xi_0)$ допускает
собственное значение с положительной вещественной частью.
Рассуждая аналогично предыдущему пункту, приходим к существованию
такого $\eps_0>0,$ что при $\eps\in(0,\eps_0]$ матрица
$u'_\xi(T,\xi_\eps,\eps)$ допускает собственное значение
$\lambda_\eps$ с б/'ольшей единицы вещественной частью.
Зафиксируем $\eps\in(0,\eps_0].$ На основании теоремы
Гробмана-Хартмана (см., напр., \cite{palis}, теорема~4.1)
существует $\lambda>0$ и локальный гомеоморфизм
$g:B_\alpha(\xi_\eps)\to\mathbb{R}^n$ такой, что
\begin{equation}\label{IZ}
  u(T,\xi,\eps)=g^{-1}(u'_\xi(T,\xi_\eps,\eps)g(\xi)),\quad\xi\in
  B_\alpha(\xi_\eps),
\end{equation}
соответственно для $p$-й степени оператора $\xi\mapsto
u(T,\cdot,\eps)$ имеем
$$
  u^p(T,\xi,\eps)=g^{-1}((u'_\xi(T,\xi_\eps,\eps))^p g(\xi)),\quad\xi\in
  B_\alpha(\xi_\eps).
$$
Пусть $l\in\mathbb{R}^n$ -- собственный вектор матрицы
$u'_\xi(T,\xi_\eps,\eps),$ соответствующий собственному значению
$\lambda_\eps$ и такой, что $g^{-1}(u'_\xi(T,\xi_\eps,\eps)l)$
определено. Из (\ref{IZ}) имеем
$g(\xi_\eps)=u'_\xi(T,\xi_\eps,\eps)g(\xi_\eps),$ то есть
$g(\xi_\eps)=0$ и, значит, $g^{-1}(l)\not=\xi_\eps.$ Поэтому, для
доказательства неустойчивости достаточно предъявить такую
сходящуюся к $\xi_\eps$ последовательность
$\{\zeta_p\}_{p\in\mathbb{N}},$ что
$$
  u^p(T,\zeta_p,\eps)=g^{-1}(l),\quad p\in\mathbb{N}.
$$
Требуемой последовательностью является, например,
$$
   \zeta_p=g^{-1}\left(\frac{1}{\lambda^p}l\right),\quad
   p\in\mathbb{N}.
$$
Действительно, так как $\xi_\eps=g^{-1}(0),$ то $\zeta_p\to
\xi_\eps$ при $p\to\infty,$ и, далее,
\begin{eqnarray*}
  u^p(T,\zeta_p,\eps)&=& g^{-1}\left((u'_\xi(T,\xi_\eps,\eps))^p
  g\left(g^{-1}\left(\frac{1}{\lambda^p}l\right)\right)\right)=
\\&=&g^{-1}\left((u'_\xi(T,\xi_\eps,\eps))^p\frac{1}{\lambda^p}l\right)=
g^{-1}\left(\lambda^p\frac{1}{\lambda^p}l\right)=g^{-1}(l).
\end{eqnarray*}
Теорема доказана полностью. \qed

Так как замена (\ref{ZAM}) $T$-периодична, то в условиях
теоремы~\ref{th1} функция $x_\eps(t)=x(t,u_\eps(t),0)$ является
$T$-периодическим решением системы (\ref{ps}) и решение $x_\eps$
устойчиво или неустойчиво вместе с $u_\eps.$ Для удобства ссылок
сформулируем это утверждение в виде теоремы.

\begin{theorem}\label{th2} Пусть выполнены условия теоремы~\ref{th1}. Тогда
утверждения 1), 2) и 3) этой теоремы имеют место и для системы
(\ref{ps}).
\end{theorem}

\section{Колебания скорости тела, перемещаемого периодическими вибрациями}

В этой секции теорема~\ref{th2} иллюстрируется на примере
доказательства периодичности и устойчивости колебаний скорости
тела в механической модели из рисунка~2.
\begin{figure}[htb]
\centerline{\includegraphics[scale=0.37]{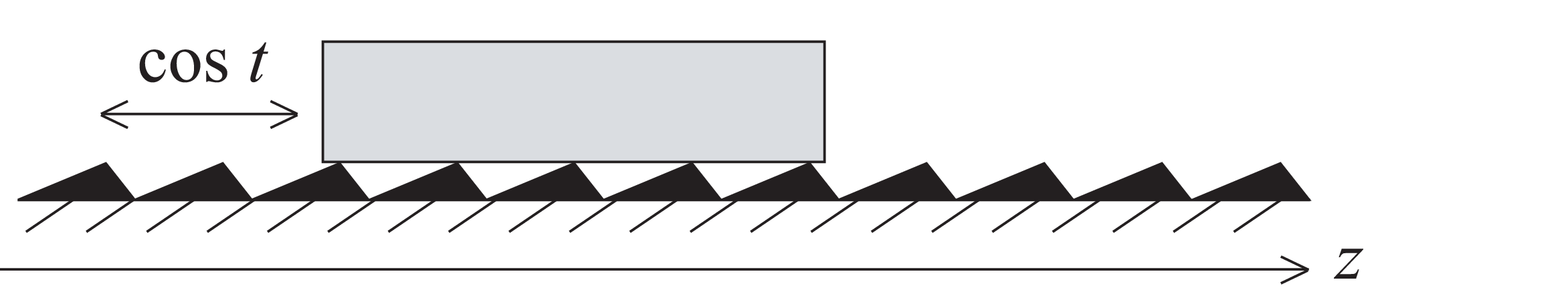}} 
\label{fidl} {\footnotesize Рис. 2. Механическая система, в
которой сила сухого трения имеет значение  $-\eps a<0$ при
движении тела вправо и значение $\eps b>0$ при движении тела
влево, где $a\not=b.$ Движение происходит за счет горизонтальной
вибрации с амплитудой $\cos t.$}
\end{figure}
Уравнение движения тела записывается (см. \cite{fidlin}) в виде
\begin{equation}\label{ex}
  \ddot z=\cos t-a\eps E(\dot z)+b\eps E(-\dot z),\ \ \mbox{где}\
E(\dot z)=({\rm sign}(\dot z)+1)/2.
\end{equation}
Замена $x=\dot z$ приводит систему (\ref{ex}) к системе вида
(\ref{ps})
\begin{equation}\label{ex1}
  \dot x=\cos t-a\eps E(x)+b\eps E(-x).
\end{equation}
Значит, $S=\{0\},$ $x(t,\xi_0,0)=-\xi_0 \sin t$ и условия
теоремы~\ref{th2} выполнены для любого
$\xi_0\in\mathbb{R}\backslash\{0\}.$  Функция $\overline{f}$ имеет
вид (см. \cite{fidlin})
$$
  \overline{f}(\xi)=-4(a+b)\arcsin(\xi)+2\pi(a-b),
$$
откуда $\xi_0=\sin\left(\frac{a-b}{a+b}\pi\right)$ и
$(\overline{f})'(\xi_0)=-2\frac{a+b}{\left|\cos\left(\frac{a-b}{a+b}\pi\right)\right|}<0.$
Следовательно, при $a>b$ ($a<b$) тело движется вправо (влево) с
$\pi$-периодически изменяющейся асимптотически устойчивой
скоростью.

\vskip0.2cm

Работа поддержана грантом BF6M10 Роснауки и CRDF (программа BRHE)
и грантом MK-1620.2008.1 Президента РФ молодым кандидатам наук.
Исследования проведены в ходе стажировки автора в Институте
Проблем Управления РАН под руководством проф.~В.Н.~Тхая и
финансируемой грантом РФФИ 08-01-90704-моб\_ст.

\end{document}